\documentclass{article}
\usepackage[utf8]{inputenc}
\usepackage{stmaryrd}
\usepackage{mathrsfs}
\usepackage[centertags]{amsmath}
\usepackage{amsfonts, dsfont}
\usepackage{amssymb}
\usepackage{amsthm}
\usepackage[OT2,T1]{fontenc}
\DeclareSymbolFont{cyrletters}{OT2}{wncyr}{m}{n}
\DeclareMathSymbol{\Sha}{\mathalpha}{cyrletters}{"58}

\input xy
\xyoption{all}

\theoremstyle{plain}
\newtheorem{thm}{Theorem}[section]
\newtheorem{cor}[thm]{Corollary}
\newtheorem{lem}[thm]{Lemma}
\newtheorem{prop}[thm]{Proposition}

\theoremstyle{definition}
\newtheorem{defn}[thm]{Definition}
\newtheorem*{remark}{Remark}
\newtheorem*{ack}{Acknowledgments}

\newcommand{\bd}{\begin{defn}}
\newcommand{\ed}{\end{defn}}
\newcommand{\bl}{\begin{lem}}
\newcommand{\el}{\end{lem}}
\newcommand{\bp}{\begin{prop}}
\newcommand{\ep}{\end{prop}}
\newcommand{\bt}{\begin{thm}}
\newcommand{\et}{\end{thm}}
\newcommand{\bc}{\begin{cor}}
\newcommand{\ec}{\end{cor}}
\newcommand{\br}{\begin{remark}}
\newcommand{\er}{\end{remark}}
\newcommand{\bdi}{\begin{diagram}}
\newcommand{\edi}{\end{diagram}}
\newcommand{\beq}{\begin{equation}}
\newcommand{\eeq}{\end{equation}}
\newcommand{\ba}{\begin{array}}
\newcommand{\ea}{\end{array}}
\newcommand{\bpf}{\begin{proof}}
\newcommand{\epf}{\end{proof}}

\newcommand{\R}{\mathds{R}}
\newcommand{\Z}{\mathds{Z}}
\newcommand{\Q}{\mathds{Q}}
\newcommand{\Zp}{\mathds{Z}_{p}}
\newcommand{\Qp}{\mathds{Q}_{p}}
\newcommand{\Op}{\mathcal{O}}
\newcommand{\al}{\alpha}
\newcommand{\be}{\beta}
\newcommand{\Ga}{\Gamma}
\newcommand{\ga}{\gamma}

\newcommand{\la}{\lambda}

\DeclareMathOperator{\Sel}{Sel} \DeclareMathOperator{\Gal}{Gal}
\DeclareMathOperator{\Hom}{Hom} \DeclareMathOperator{\rank}{rank}
\DeclareMathOperator{\corank}{corank}

\newcommand{\cts}{\mathrm{cts}}
\newcommand{\cyc}{\mathrm{cyc}}
\newcommand{\M}{\mathfrak{M}}

\newcommand{\ot}{\otimes}
\newcommand{\ilim}{\displaystyle \mathop{\varinjlim}\limits}
\newcommand{\plim}{\displaystyle \mathop{\varprojlim}\limits}

\newcommand{\coker}{\mathrm{coker}\,}

\newcommand{\lra}{\longrightarrow}

\newcommand{\ps}[1]{\llbracket #1 \rrbracket}

\begin{document}
\title{Some remarks on Kida's formula when $\mu\neq 0$}
 \author{
  Meng Fai Lim\footnote{School of Mathematics and Statistics $\&$ Hubei Key Laboratory of Mathematical Sciences,
Central China Normal University, Wuhan, 430079, P.R.China.
 E-mail: \texttt{limmf@mail.ccnu.edu.cn}} }
\date{}
\maketitle

\begin{abstract} \footnotesize
\noindent
The Kida's formula in classical Iwasawa theory relates the Iwasawa $\la$-invariants of $p$-extensions
of number fields. Analogue of this formula was subsequently established for the Iwasawa $\la$-invariants of Selmer groups under an appropriate $\mu=0$ assumption. In this paper, we give a conceptual (but conjectural) explanation that such a formula should also hold when $\mu\neq 0$. The conjectural component comes from the so-called $\M_H(G)$-conjecture in noncommutative Iwasawa theory.


\end{abstract}

\section{Introduction}

The classical Kida's formula relates the $\la$-invariants of CM-extensions of cyclotomic $\Zp$-extensions of number field (see \cite{Ki}; also see \cite{Iw}). An analogue of this formula for Selmer groups of elliptic curves was later discovered by Hachimori-Matsuno \cite{HM}. Subsequently, this has been generalized to Selmer groups of ordinary modular forms and more general Galois representations (for instances, see \cite{CDLSS,PW06}). A common phenomenon in deriving such a formula is that one has to assume an appropriate $\mu=0$ assumption. It is then natural to ask if such a formula still holds without the $\mu=0$ assumption. Indeed if the Galois representation possess a Galois-invariant lattice whose Selmer group satisfies $\mu=0$ hypothesis, one does have an analogous formula (see \cite[Remark 3.2.3]{Gr11}). However, this observation cannot be carried out for more general situation. Indeed, the work of Drinen \cite{Dr} (also see \cite{DD15}) has exhibited examples of Galois representations whose Galois-invariant lattices do not satisfy the $\mu=0$ assumption.

The aim of this paper is to at least provide a conceptual (but conjectural) explanation to the validity of Kida's formula when $\mu\neq 0$.
We briefly explain the idea behind our approach here, leaving the details to the body of the paper. Let $F_{\infty}$ be a strongly admissible pro-$p$ extension of $F$, where we write $H=\Gal(F_{\infty}/F^{\cyc})$ and $\Ga=\Gal(F^{\cyc}/F)$. Suppose for now we are in the $\mu=0$ situation. In other words, our dual Selmer group $X(A/F^{\cyc})$ is finitely generated over $\Zp$ (for simplicity of exposition here, we assume our Galois module has coefficients in $\Zp$; see the body of the paper for a more general context). A by-now standard argument then shows that the dual Selmer group $X(A/F_{\infty})$ is finitely generated over $\Zp\ps{H}$ (for instances, see \cite[Theorem 6.4(ii)]{CH} or \cite[Theorem 2.1]{CS12}). Furthermore, there is a well-known formula relating the $\Zp\ps{H}$-rank of $X(A/F_{\infty})$ with $\Zp$-rank of $X(A/F^{\cyc})$ (see \cite[Theorem 5.4]{HS} or \cite[Corollary 2.12]{Ho}). There are two approaches towards proving this formula. The first approach is via a descent argument comparing the cyclotomic Selmer groups under a finite $p$-base change. This approach made heavy usage of Kida's formula (see \cite[Theorem 3.1]{HV} or \cite[Theorem 5.4]{HS}). A second approach is a homological descent argument due to Howson (cf. \cite[Corollary 2.12]{Ho}). We emphasize that Howson's approach avoids use of Kida's formula. An important observation, which does not seem to be written down in literature (although the author believes this might be known) is that the above said rank formula can be used to recover Kida's formula. We briefly say a few words on this observation. For a subgroup $H'$ of $H$ with finite index, one has an equality $|H:H'|\rank_{\Zp\ps{H}}\big(X(A/F_{\infty})\big) = \rank_{\Zp\ps{H'}}\big(X(A/F_{\infty})\big)$. As seen above, the quantity $\rank_{\Zp\ps{H}}\big(X(A/F_{\infty})\big)$ is related to $\Zp$-rank of $X(A/F^{\cyc})$, and on the other hand, $\rank_{\Zp\ps{H'}}\big(X(A/F_{\infty})\big)$ can be related to $\Zp$-rank of $X(A/L^{\cyc})$, where $L$ is a finite $p$-extension of $F$ contained in $F_{\infty}$ with $L^{\cyc} = (F_{\infty})^{H'}$. Combining these three relations, we obtain Kida's formula.

Now we turn to the $\mu\neq 0$ situation. Inspired by the above discussion, the natural first step towards establishing a Kida formula is then to note that under the validity of the $\M_H(G)$-conjecture, Howson's proof can be adapted to establish an analogous formula relating the $\Zp\ps{H}$-rank of $X(A/F_{\infty})/X(A/F_{\infty})(p)$ over $F_{\infty}$ with $\Zp$-rank of $X(A/F^{\cyc})/X(A/F^{\cyc})(p)$ (see Proposition \ref{rank formula}). Consequently, we may combine this formula with the above strategy to obtain Kida's formula (see Theorem \ref{Kida formula theorem}).
Although this gives a conceptual explanation to Kida's formula, we have to confess that we are not able to establish the validity of such a formula (much less the $\M_H(G)$-conjecture) unconditionally.
Nevertheless, we at least content ourselves with recovering the original Kida formula for the $\mu=0$ case via the approach given here (see Theorem \ref{Kida formula theorem mu=0} and also discussion in Section \ref{elliptic review}).

We now give a brief description of the contents of each section of the paper. Section \ref{algebra} introduces some algebraic notations and results required for the paper. In Section \ref{Selmer}, we introduce the (strict) Selmer group following Greenberg \cite{G89} and prove a formula relating the ranks for the dual Selmer groups of the admissible extension and cyclotomic extension. Section \ref{Kida section} is where we prove our main result. Finally, in Section \ref{elliptic review}, we revisit the elliptic curve situation as illustration of our result.

\begin{ack}
The author like to thank John Coates for his interest and comments on the paper.
He would also like to thank the anonymous referee for many useful comments and suggestions on the article.  Some part of the research of this article was conducted when the author was visiting the National University of Singapore and the National Center for Theoretical Sciences in Taiwan, and he would like to acknowledge the hospitality and conducive working conditions provided by these institutes.
The author's research is supported by the National Natural Science Foundation of China under Grant No. 11550110172 and Grant No. 11771164. \end{ack}

\section{Algebraic Preliminaries} \label{algebra}

Throughout the paper, $p$ will denote a fixed prime. Let $\Op$ be the ring of integers of a finite extension of $\Qp$. Fix a local parameter
$\pi$ for $\Op$ and denote by $k$ the residue field of $\Op$. Let $G$ be a compact pro-$p$ $p$-adic
Lie group without $p$-torsion. The completed group algebra of $G$ over $\Op$ is given by
 \[ \Op\ps{G} = \plim_U \Op[G/U], \]
where $U$ runs over the open normal subgroups of $G$ and the inverse
limit is taken with respect to the canonical projection maps. It is well known that $\Op\ps{G}$ is
an Auslander regular ring (cf. \cite[Theorem 3.26]{V02} or \cite[Theorem A.1]{LimFine}) with no zero divisors (cf.\
\cite{Neu}). Therefore, it admits a skew field $Q(G)$ which is flat
over $\Zp\ps{G}$ (see \cite[Chapters 6 and 10]{GW} or \cite[Chapter
4, \S 9 and \S 10]{Lam}). Thanks to this, we can define the notion of $\Op\ps{G}$-rank of a finitely generated $\Op\ps{G}$-module $M$, which is in turn given by
$$ \rank_{\Op\ps{G}}(M)  = \dim_{Q(G)} (Q(G)\ot_{\Op\ps{G}}M). $$
The module $M$ is then said to be a
\textit{torsion} $\Op\ps{G}$-module if $\rank_{\Op\ps{G}} (M) = 0$.

The completed group algebra $k\ps{G}$ is defined similarly as above.
By \cite[Theorem 3.30(ii)]{V02}, the ring $k\ps{G}$ is Auslander
regular. Furthermore, it has no zero divisors (cf. \cite[Theorem
C]{AB}). Therefore, one has a notion of $k\ps{G}$-rank defined as
above. Similarly, we shall say that
the module $N$ is a \textit{torsion} $k\ps{G}$-module if
$\rank_{k\ps{G}}(N) = 0$.

For a given finitely generated $\Op\ps{G}$-module $M$, denote by
$M(\pi)$ the $\Op\ps{G}$-submodule of $M$ consisting of elements
of $M$ which are annihilated by some power of $\pi$. Since the ring
$\Op\ps{G}$ is Noetherian, the module $M(\pi)$ is also finitely generated
over $\Op\ps{G}$. Thus, there exists a sufficiently large integer $r$ so
that $\pi^r$ annihilates $M(\pi)$. Following \cite[Formula (33)]{Ho}, we
define
  \[\mu_{\Op\ps{G}}(M) = \sum_{i\geq 0}\rank_{k\ps{G}}\big(\pi^i
   M(\pi)/\pi^{i+1}\big). \]
(For another alternative, but equivalent, definition, see
\cite[Definition 3.32]{V02}.) Note that the sum on
the right is a finite one by the above discussion.

Now if $M$ is a finitely generated $\Op\ps{G}$-module, then its homology groups $H_i(G,M)$ are finitely generated over $\Op$ (see \cite[Proof of Theorem 1.1]{Ho} or \cite[Lemma 3.2.3]{LS}). Hence the quantity $\rank_{\Op}\big(H_i(G,M)\big)$ is well-defined. In view of this observation, we can now state the following result of Howson  (see \cite[Theorem 1.1]{Ho} or \cite[Lemma 4.3]{LimFine}).

\bp[Howson] \label{Howson}
Let $M$ be a finitely generated $\Op\ps{G}$-module. Then we have
\[ \rank_{\Op\ps{G}}(M) = \sum_{i=0}^d(-1)^i\rank_{\Op}\big(H_i(G,M)\big),\]
where here $d$ denotes the dimension of the $p$-adic group $G$.
\ep

In this paper, we are mostly interested in the class of $p$-adic Lie groups $G$ which contains a closed normal subgroup $H$ such that $G/H\cong\Zp$. The following relative formula for the $\mu$-invariant will be of use later.

\bl \label{mu compare Mf}
Let $G$ be a pro-$p$ compact $p$-adic group without $p$-torsion which contains a closed normal subgroup $H$ with the property that $\Ga:=G/H\cong \Zp$. Let $M$ be a finitely generated $\Op\ps{G}$-module. Then we have
\[ \mu_{\Op\ps{G}}(M) = \sum_{i=0}^d(-1)^i\mu_{\Op\ps{\Ga}}\big(H_i(H,M(\pi))\big).\]
\el

\bpf
See \cite[Lemma 2.2]{LimMHG}.
\epf

For a group $G$ defined as above, a finitely generated $\Op\ps{G}$-module $M$ is said to
\textit{satisfy the $\M_H(G)$-property} if $M_f:=M/M(\pi)$ is
finitely generated over $\Op\ps{H}$. It has been conjectured for
certain Galois representations coming from abelian varieties with
good ordinary reduction at $p$ or cuspidal eigenforms with good
ordinary reduction at $p$, the dual Selmer group associated to such
a Galois representation satisfies the $\M_H(G)$-property (see
\cite{BV, CFKSV, CS12, LimMHG}).

For the subsequent discussion, we shall write $M_H$ for the largest quotient of $M$ on which $H$ acts trivially. Note that $M_H = H_0(H,M)$. We now record another useful lemma.

\bl \label{rank compare Mf}
Let $G$ be a pro-$p$ compact $p$-adic group without $p$-torsion which contains a closed normal subgroup $H$ such that $\Ga:=G/H\cong \Zp$. Let $M$ be a finitely generated $\Op\ps{G}$-module which satisfies the $\M_H(G)$-property and has the property that $H_i(H,M)$ is finitely generated over $\Op$ for all $i\geq 1$. Then we have
\[ \mu_{\Op\ps{G}}(M) = \mu_{\Op\ps{\Ga}}\big(M_H\big)\]
and \[ \rank_{\Op}\big(H_i(H, M)\big) =\rank_{\Op}\big(H_i(H, M_f)\big)\]
for every $i\geq 1$.
\el

\bpf
 Taking $H$-homology
of the following short exact sequence
\[ 0\lra M(\pi)\lra M\lra M_f\lra 0,\]
yields an exact sequence
\[   H_{i+1}(H, M_f)\lra H_{i}(H, M(\pi))\lra H_{i}(H,
M) \stackrel{f_i}\lra H_i(H, M_f)\lra H_{i-1}(H, M(\pi))\]
for $i\geq 1$.  As $M_f$ is finitely generated over $\Op\ps{H}$, the group $H_i(H, M_f)$ is therefore finitely generated over $\Op$ for every $i$ (cf. \cite[Proof of Theorem 1.1]{Ho} or \cite[Lemma 3.2.3]{LS}). Combining this with the hypothesis that $H_i(H,M)$ is finitely generated over $\Op$ for all $i\geq 1$, we see that $H_i(H,M(\pi))$ is finitely generated over $\Op$ for $i\geq 1$. As noted before, there exists a sufficiently large $r$ such that $\pi^r$ annihilates $M(\pi)$, and hence all the $H$-homology groups of $M(\pi)$. Therefore, $H_i(H,M(\pi))$ must be finite for $i\geq 1$. From this and the above exact sequence, we have the second equality of the lemma. Furthermore, the finiteness implies $\mu_{\Op\ps{\Ga}}\big(H_i(H,M(\pi))\big) =0$ for $i\geq 1$. Putting these into Lemma \ref{mu compare Mf}, we obtain
\[ \mu_{\Op\ps{G}}(M) = \mu_{\Op\ps{\Ga}}\big(M(\pi)_H\big).\]
On the other hand, from the exact sequence
\[   H_{1}(H, M_f)\lra  M(\pi)_H\lra M_H \lra H_0(H, M_f)\lra 0,\]
and the observation that each $H_i(H,M_f)$, being finitely generated over $\Op$, has trivial $\mu_{\Op\ps{\Ga}}$-invariants, we see that
\[ \mu_{\Op\ps{\Ga}}\big(M(\pi)_H\big) = \mu_{\Op\ps{\Ga}}\big(M_H\big).\]
Combining this with the above, we obtain the first equality of the lemma.

Finally, the finiteness of $H_i(H,M(\pi))$ (for $i\geq 1$) immediately gives the second equality for $i\geq 2$. For $i=1$, it remains to verify that $\mathrm{coker} f_1$ to be finite. But this follows from the facts that $\mathrm{coker} f_1$ is finitely generated over $\Op$ (being a quotient of $H_1(H, M_f)$) and is annihilated by $\pi^r$ (being a submodule of $M(\pi)_H$). The proof of the lemma is now completed.
\epf

\section{Selmer groups} \label{Selmer}

Throughout this section, we let $F$ be a number
field which is further assumed to have no real primes if $p=2$. As before, denote by $\Op$ the ring of integers of a fixed finite extension of $\Qp$. In this section, we define the strict Selmer groups associated to certain datum in the sense of Greenberg \cite{G89}. As a start, we introduce the axiomatic conditions on our datum which is denoted by $\big(A, \{A_v\}_{v|p}\big)$ and satisfies all of the following four conditions \textbf{(C1)-(C4)}.

\begin{enumerate}
 \item[(\textbf{C1})] $A$ is a
cofinitely generated cofree $\Op$-module of $\Op$-corank $d$ with a
continuous, $\Op$-linear $\Gal(\bar{F}/F)$-action which is
unramified outside a finite set of primes of $F$.

 \item[(\textbf{C2})] For each prime $v$ of $F$ above $p$, there is a
distinguished $\Gal(\bar{F}_v/F_v)$-submodule $A_v$ of $A$ which is cofree of
$\Op$-corank $d_v$.

 \item[(\textbf{C3})] For each real prime $v$ of $F$, $A_v^+:=
A^{\Gal(\bar{F}_v/F_v)}$  is cofree of
$\Op$-corank $d^+_v$.

\item[(\textbf{C4})] The following equality
  \[\sum_{v|p} (d-d_v)[F_v:\Qp] = dr_2(F) +
 \sum_{v~\mathrm{real}}(d-d^+_v)\]
holds. Here $r_2(F)$ denotes the number of complex primes of $F$.
\end{enumerate}

As we need to work with Selmer groups defined over a tower of number fields, we need to consider the base change of our datum which we now do. For a finite extension $L$ of $F$, the base change of our datum $\big(A,
\{A_w\}_{w|p} \big)$ over $L$ is given as follows:

(1) $A$ can be viewed as a $\Gal(\bar{F}/L)$-module via restriction of the Galois action.

(2) For each prime $w$
of $L$ above $p$, we set $A_w =A_v$, where $v$ is the prime of $F$
below $w$, and view it as a $\Gal(\bar{F}_v/L_w)$-module via the appropriate restriction.
Note that we then have $d_w= d_v$.

(3) For each real prime $w$ of $L$ which lies above a real prime $v$ of $F$, we set
$A_w^+= A^{\Gal(\bar{F}_v/F_v)}$ and write $d^+_w
= d^+_v$.

We now record the following lemma which gives some
sufficient conditions for equality in \textbf{(C4)} to hold for
the datum $\big(A, \{A_w\}_{w|p}\big)$ over $L$.

\bl \label{data base change} Suppose that $\big(A, \{A_v\}_{v|p}\big)$ is a datum defined over $F$.
 Suppose further that at least one of the following statements holds.
 \begin{enumerate}
\item[$(i)$] $[L:F]$ is odd.

\item[$(ii)$] $F$ has no real primes.
 \end{enumerate}
Then the datum $\big(A, \{A_w\}_{w|p}\big)$ obtained by base change satisfies \textbf{(C1)-(C4)}. In particular, we have the equality
 \[ \sum_{w|p} (d-d_w)[L_w:\Qp] = dr_2(L) +
 \sum_{w~\mathrm{real}}(d-d^+_w).\]
 \el

\bpf See \cite[Lemma 3.0.1]{LimCMu}. \epf

We shall impose one more condition on our datum.

\medskip\textbf{(Fin)} $H^0(L^{\cyc}, A^*)$ is finite for every $L$ contained in $F_{\infty}$. Here $A^* =
\Hom_{\cts}(T_{\pi}(A),\mu_{p^{\infty}})$, where $T_{\pi}(A) =
\plim_i A[\pi^i]$.

\medskip
We mention two basic examples of our datum.

(i) $A = \mathcal{A}[p^{\infty}]$, where $\mathcal{A}$ is an abelian
variety defined over an arbitrary finite extension $F$ of $\Q$ with
good ordinary reduction at all places $v$ of $F$ dividing $p$. For
each $v|p$, it follows from \cite[P. 150-151]{CG} that we have a
$\Gal(\bar{F}_v/F_v)$-submodule $A_v$ characterized by
the property that $A/A_v$ is the maximal
$\Gal(\bar{F}_v/F_v)$-quotient of $\mathcal{A}[p^{\infty}]$ on which
some subgroup of finite index in the inertia group $I_v$ acts
trivially. It is not difficult to verify that \textbf{(C1)-(C4)} are
satisfied. The condition \textbf{(Fin)} is a well-known
consequence of a theorem of Imai \cite{Imai}.

(ii) Let $V$ be the Galois representation attached to a primitive
Hecke eigenform $f$ for $GL_2 /\Q$, which is ordinary at $p$,
relative to some fixed embedding of the algebraic closure of $\Q$
into $\Qp$. By the work of Mazur-Wiles \cite{MW}, $V$ contains a
one-dimensional $\Qp$-subspace $V_v$ invariant under
$\Gal(\bar{\Q}_p/\Qp)$ with the property that the inertial subgroup
$I_p$ acts via a power of the cyclotomic character on $V_v$ and
trivially on $V/V_v$. By compactness, $V$ will always contain a free
$\Op$-submodule $T$, which is stable under the action of
$\Gal(\bar{F} /F)$. For such an $\Op$-lattice $T$, we write $A =
V/T$ and $A_v = V_v/ (T\cap V_v)$. The condition \textbf{(Fin)} is shown in
the proof of \cite[Lemma 2.2]{Su}.

We return to general discussion.
A Galois extension $F_{\infty}$ of $F$ is said to be a strongly admissible pro-$p$ $p$-adic Lie
extension of $F$ if (i) $\Gal(F_{\infty}/F)$ is a compact pro-$p$ $p$-adic Lie group with no $p$-torsion, (ii) $F_{\infty}$ contains the cyclotomic $\Zp$-extension
$F^{\cyc}$ of $F$ and (iii) $F_{\infty}$ is unramified outside a finite set of primes. We shall always write $G = \Gal(F_{\infty}/F)$, $H = \Gal(F_{\infty}/F^{\cyc})$ and $\Ga =\Gal(F^{\cyc}/F)$. By Lemma \ref{data base change}, the datum $\big(A, \{A_w\}_{w|p}, \{A^+_w\}_{w|\R} \big)$ obtained by base changing to any intermediate subextension of $F_{\infty}/F$ satisfies \textbf{(C1)-(C4)}. (Recall that when $p=2$, our standing assumption is that $F$ has no real primes.)

We can now define the strict Selmer group associated to our datum following Greenberg \cite{G89}. Let $S$ be a finite set of
primes of $F$ which contains all the primes above $p$, the ramified
primes of $A$, the ramified
primes of $F_{\infty}/F$ and all the infinite primes of $F$. Denote by $F_S$ the
maximal algebraic extension of $F$ unramified outside $S$ and write
$G_S(\mathcal{L}) = \Gal(F_S/\mathcal{L})$ for every algebraic
extension $\mathcal{L}$ of $F$ which is contained in $F_S$. Let $L$
be a finite extension of $F$ contained in $F_{\infty}$. For a prime $w$ of $L$ lying over $S$, set
\[ H^1_{str}(L_w, A)=
\begin{cases} \ker\big(H^1(L_w, A)\lra H^1(L_w, A/A_w)\big), & \text{\mbox{if} $w$
 divides $p$},\\
 \ker\big(H^1(L_w, A)\lra H^1(L^{ur}_w, A)\big), & \text{\mbox{if} $w$ does not divide $p$,}
\end{cases} \]
 where $L_w^{ur}$ is the maximal unramified extension of $L_w$.
 We then define
 \[ J_v(A/L) = \bigoplus_{w|v}H^1_s(L_w, A),\]
where $w$ runs over the (finite) set of primes of $L$ above $v$ and $H^1_s(L_w, A) = H^1(L_w, A)/H^1_{str}(L_w, A)$. If
$\mathcal{L}$ is an infinite extension of $F$, we define
\[ J_v(A/\mathcal{L}) = \ilim_L J_v(A/L),\]
where the direct limit is taken over all finite extensions $L$ of
$F$ contained in $\mathcal{L}$. For a (possibly infinite) extension $K$ of $F$ contained in $F_{\infty}$, the (strict) Selmer group is defined by
\[ S^{str}(A/K) := \Sel^{str}(A/K) := \ker\left( H^1(G_S(K),A)\lra
\bigoplus_{v\in S}J_v(A/K)\right).\]
It follows from the definition that $S^{str}(A/K) = \ilim_L S^{str}(A/L)$, where the direct limit is taken
with respect to the natural maps, as $L$ varies over all finite subextensions
of the base field $F$ in the larger extension $K$. Write $X(A/K)$ for its Pontryagin dual and  $X_f(A/K) =X(A/K)/X(A/K)(\pi)$.

Suppose that $X_f(A/F_{\infty})$ satisfies the $\M_H(G)$-property. It then makes sense to speak of the quantity $\rank_{\Op\ps{H}}\big(X_f(A/F_{\infty})\big)$. Also, it follows from \cite[Proposition 2.5]{CS12} that $X(A/F^{\cyc})$ is torsion over $\Op\ps{\Ga}$ and so the quantity $\rank_{\Op}(X_f(A/F^{\cyc}))$ is well-defined. We can therefore state the following.

\bp \label{rank formula}
Let $F_{\infty}$ be a strongly admissible pro-$p$ Lie extension of $F$. Suppose that the data $\big(A, \{A_v\}_{v|p}, \{A^+_v\}_{v|\R} \big)$ satisfies $\textbf{(C1)-(C4)}$ and $\textbf{(Fin)}$. Assume that $X(A/F_{\infty})$ satisfies the $\M_H(G)$-property. Then we have
\[ \mu_{\Op\ps{G}}\big(X(A/F_{\infty})\big) = \mu_{\Op\ps{\Ga}}\big(X(A/F^{\cyc})\big)\]
and
\[ \rank_{\Op\ps{H}}\big(X_f(A/F_{\infty})\big) = \rank_{\Op}\big(X_f(A/F^{\cyc})\big) - \corank_{\Op}\big(A(F^{\cyc})\big) +
\sum_{\substack{w\in S(F^{\cyc}),
\\ \dim H_w\geq 1}}\corank_{\Op}\big(D_w(F^{\cyc}_w)\big).\]
Here $D_w$ denotes $A/A_w$ or $A$ accordingly as $w$ divides $p$ or not.
\ep

Proposition \ref{rank formula} has been proved for an elliptic curve $E$ under the stronger assumption that $X(E/F_{\infty})$ is finitely generated over $\Zp\ps{H}$ (see
\cite[Corollary 6.10]{CH}, \cite[Theorem 5.4]{HS}, \cite[Theorem
3.1]{HV} and \cite[Theorem 2.8]{Ho}). The approach for the proof in \cite{CH, HS, HV} builds on Kida's formula which we want to avoid for our purposes. On the other hand, the proof given in \cite[Theorem 2.8]{Ho} does not assume the validity of Kida's formula. The goal of Proposition \ref{rank formula} is to show that Howson's approach works under the $\M_H(G)$-property which the remainder of the section is devoted to.

As a start, we have the following lemma.

\bl \label{short exact sequences}
Retaining the assumptions of Proposition \ref{rank formula}, we have short exact sequences
\[ 0 \lra S^{str}(A/F^{\cyc})\lra H^1(G_S(F^{\cyc}), A)  \lra \bigoplus_{v\in S} J_v(A/F^{\cyc})\lra 0\]
and
\[ 0 \lra S^{str}(A/F_{\infty})\lra H^1(G_S(F_{\infty}), A)  \lra \bigoplus_{v\in S} J_v(A/F_{\infty})\lra 0.\]
\el

\bpf
Since $X(A/F_{\infty})$ satisfies the $\M_H(G)$-property, it follows from \cite[Proposition 2.5]{CS12} that for every finite extension $L$ of $F$ contained in $F_{\infty}$, $X(A/L^{\cyc})$ is torsion over $\Op\ps{\Ga_L}$, where $\Ga_L=\Gal(L^{\cyc}/L)$. In view of $\textbf{(Fin)}$, we may apply a similar argument to that in \cite[Proposition 3.3]{LimMHG} to obtain a short exact sequence
\[ 0 \lra S^{str}(A/L^{\cyc})\lra H^1(G_S(L^{\cyc}), A)     \lra \bigoplus_{v\in S} J_v(A/L^{\cyc})\lra 0.\]
In particular, this yields the first short exact sequence by taking $L=F$.
On the other hand, by taking direct limit over $L$, we obtain the second short exact sequence.
\epf

The next two lemmas are concerned with the $H$-homology of global cohomology groups and local cohomology groups.

\bl \label{global calculation}
Retain the assumptions of Proposition \ref{rank formula}. Then $H^i\big(H, H^1\big(G_S(F_{\infty}),A)\big)$ is cofinitely generated over
$\Op$ for every $i\geq 1$. Moreover, we have an exact sequence
\[ 0 \lra H^1(H,
A(F_{\infty}))\lra H^1(G_S(F^{\cyc}), A) \lra H^1(G_S(F_{\infty}), A)^H\lra H^2(H,
A(F_{\infty}))\lra 0\]
and isomorphisms \[
H^i\big(H, H^1(G_S(F_{\infty}), A)\big)\cong H^{i+2}\big(H,
A(F_{\infty})\big) \mbox{ for } i\geq 1.\]
\el

\bpf
Since $X(A/F_{\infty})$ is assumed to satisfy the $\M_H(G)$-property, it follows from \cite[Proposition 2.5]{CS12} that for every finite extension $L$ of $F$ contained in $F_{\infty}$, $X(A/L^{\cyc})$ is torsion over $\Op\ps{\Ga_L}$, where $\Ga_L=\Gal(L^{\cyc}/L)$. Via similar arguments to those in \cite[Proposition 3.3 and Corollary 3.4]{LimMHG}, we have that
$H^2(G_S(F^{\cyc}),A)=0$ and $H^2(G_S(F_{\infty}),A)=0$. Hence the spectral sequence
\[ H^i\big(H, H^j(G_S(F_{\infty}), A)\big)\Longrightarrow
H^{i+j}(G_S(F^{\cyc}), A)\] degenerates to yield an exact sequence
\[ 0 \lra H^1(H,
A(F_{\infty}))\lra H^1(G_S(F^{\cyc}), A)     \lra H^1(G_S(F_{\infty}), A)^H\lra H^2(H,
A(F_{\infty}))\lra 0\]
and isomorphisms \[
H^i\big(H, H^1(G_S(F_{\infty}), A)\big)\cong H^{i+2}\big(H,
A(F_{\infty})\big) \mbox{ for } i\geq 1.\]
Finally, the $\Op$-cofinitely generation of the latter
groups follows from the fact that the cohomology groups $H^i(H,W)$ are cofinitely generated
$\Op$-modules for any $p$-adic Lie group $H$ and any $\Op$-cofinitely generated
$H$-module $W$. Consequently,  $H^i\big(H, H^1\big(G_S(F_{\infty}),A)\big)$ is cofinitely generated over
$\Op$ for every $i\geq 1$.
\epf

\bl \label{local calculation}
Retain the assumption of Proposition \ref{rank formula}. Then  $H^i(H, \bigoplus_{v\in S}J_v(A/F_{\infty}))$ is cofinitely generated over
$\Op$ for every $i\geq 1$. Moreover, we have an exact sequence
\[ 0 \lra \bigoplus_{w\in S(F^{\cyc})}H^1\big(H_w,
D_v(F_{\infty,w})\big)\lra \bigoplus_{v\in S}J_v(A/F^{\cyc})   \lra \left(
\bigoplus_{v\in S}J_v(A/F_{\infty})\right)^H \] \[\lra  \bigoplus_{w\in S(F^{\cyc})}H^2\big(H_w,
D_v(F_{\infty,w})\big)\lra 0\]
and isomorphisms
\[
H^i\left(H, \bigoplus_{v\in S}J_v(A/F_{\infty})\right)\cong\bigoplus_{w\in S(F^{\cyc})}H^{i+2}\big(H_w,
D_v(F_{\infty,w})\big) \mbox{ for } i\geq 1.\]
Here $D_v$ denotes $A/A_v$ or $A$ accordingly as $v$ divides $p$ or not.
\el

\bpf
This is a local version of Lemma \ref{global calculation} with a similar proof noting that
$H^2(F^{\cyc}_w,A)=0$ and $H^2(F_{\infty,w},A)=0$ by \cite[Theorem 7.1.8(i)]{NSW}.
\epf

We can now give the proof of Proposition \ref{rank formula}.

\bpf[Proof of Proposition \ref{rank formula}]
 By Lemma \ref{short exact sequences}, we have the following commutative diagram
\[  \entrymodifiers={!! <0pt, .8ex>+} \SelectTips{eu}{}\xymatrix{
    0 \ar[r]^{} & S^{str}(A/F^{\cyc}) \ar[d]_{\al} \ar[r] &  H^1(G_S(F^{\cyc}), A)
    \ar[d]_{\be}
    \ar[r] & \displaystyle\bigoplus_{v\in S}J_v(A/F^{\cyc}) \ar[d]_{\ga} \ar[r] & 0 &\\
    0 \ar[r]^{} & S^{str}(A/F_{\infty})^H \ar[r]^{} & H^1(G_S(F_{\infty}), A)^H \ar[r] & \
    \displaystyle\bigoplus_{v\in S}J_v(A/F_{\infty})^H \ar[r] &  H^1\big(H, S^{str}(A/F_{\infty})\big) \ar[r] & \cdots } \]
with exact rows.  To simplify notation, we write $W_{\infty}=
H^1(G_S(F_{\infty}), A)$ and
 $J_{\infty} = \displaystyle\bigoplus_{v\in S}J_v(A/F_{\infty})$. From the commutative diagram, we have a long exact
 sequence
 \[ \ba{c} 0\lra \ker\al \lra \ker \be
 \lra \ker \ga
 \lra \coker \al \lra \coker \be \\
  \lra \coker\ga \lra H^1\big(H, S^{str}(A/F_{\infty})\big)
 \lra H^1(H, W_{\infty})
 \lra H^1(H, J_{\infty})\lra \cdots \\
 \cdots\lra H^{i-1}(H, J_{\infty}) \lra H^i\big(H, S^{str}(A/F_{\infty})\big)
 \lra H^i(H, W_{\infty})
 \lra H^i(H, J_{\infty})\lra \cdots .\ea \]

By Lemmas \ref{global calculation} and \ref{local calculation}, the groups $\ker \be$, $\ker \ga$, $\coker \be$, $\coker\ga$, $ H^i(H, W_{\infty})$ (for $i\geq 1$) and $H^{i}(H, J_{\infty})$ (for $i\geq 1$) are cofinitely generated over $\Op$. Thus, combining this observation with the above exact sequence, we have that $\ker \al$, $\coker \al$ and $H^1\big(H, S(A/F_{\infty})\big)$ (for $i\geq 1)$ are cofinitely generated over $\Op$. Therefore, we have
\[ \mu_{\Op\ps{G}}\big(X(A/F_{\infty})\big) = \mu_{\Op\ps{\Ga}}\big(X(A/F_{\infty})_H\big) =  \mu_{\Op\ps{\Ga}}\big(X(A/F^{\cyc})\big),\]
where the first equality follows from Lemma \ref{rank compare Mf} and the second follows from the facts that $\ker \al$ and $\coker \al$ are cofinitely generated over $\Op$. This establishes the first equality of the proposition. Moreover, we have
\[ \ba{c} \corank_{\Op}(\ker\al) - \corank_{\Zp}(\coker\al)=  -\displaystyle\sum_{i\geq 1}(-1)^i\corank_{\Op}H^i(H,S^{str}(A/F_{\infty}))\hspace{2in} \\
   \hspace{1in}+\displaystyle\sum_{i\geq 1}(-1)^{i}\corank_{\Op}H^i(H,
A(F_{\infty}))-\displaystyle\sum_{\substack{w\in S(F^{\cyc}),
\\ \dim H_w\geq 1}}\left(\sum_{i\geq 1}(-1)^{i}\corank_{\Op}H^i(H_w,
D_v(F_{\infty,w}))\right), \ea \]
where here $D_v$ denotes $A/A_v$ or $A$ accordingly as $v$ divides $p$ or not.
Applying Proposition \ref{Howson} and Lemma \ref{rank compare Mf}, the right hand side is just
\[ -\displaystyle\sum_{i\geq 1}(-1)^i\rank_{\Op}H_i(H,X_f(A/F_{\infty}))-\corank_{\Op}H^0(H,
A(F_{\infty}))
 +\displaystyle\sum_{\substack{w\in S(F^{\cyc}),
\\ \dim H_w\geq 1}}\corank_{\Op}H^0(H_w,
D_v(F_{\infty,w})).
\]
Now consider the following commutative diagram
\[  \entrymodifiers={!! <0pt, .8ex>+} \SelectTips{eu}{}\xymatrix{
     & X(A/F_{\infty})(p)_H \ar[d]_{h'} \ar[r] &  X(A/F_{\infty})_H
    \ar[d]_{\al^{\vee}}
    \ar[r] & X_f(A/F_{\infty})_H \ar[d]_{h''} \ar[r]& 0 \\
    0 \ar[r]^{} & X(A/F^{\cyc})(p) \ar[r]^{} &  X(A/F^{\cyc}) \ar[r] & \
     X_f(A/F^{\cyc}) \ar[r] &  0 } \]
with exact rows. From this, we have a long exact
 sequence
 \[   \ker h' \lra \ker (\al^{\vee})
 \stackrel{f}{\lra} \ker h''
 \lra \coker h' \lra \coker (\al^{\vee}) \lra \ker h''\lra 0.  \]
Since $X(A/F_{\infty})$ satisfies $\M_H(G)$-property, it follows from Nakayama's lemma that $X_f(A/F_{\infty})_H$ is finitely generated over $\Op$. But $X_f(A/F^{\cyc})$ is also finitely generated over $\Op$ by \cite[Proposition 2.5]{CS12}. Hence $\ker  h''$ and $\coker h''$ are finitely generated over $\Op$, and we have
\[ \rank_{\Op} (\ker h'') - \rank_{\Op}(\coker h'')= \rank_{\Op}\big(X_f(A/F_{\infty})_H\big) -\rank_{\Op}(X_f(A/F^{\cyc}))\]

On the other hand, as already seen above, $\ker (\al^{\vee})$ and $\coker (\al^{\vee})$ are finitely generated over $\Op$. Hence so are $\ker f$ and $\coker h'$. But since these latter groups are $\pi$-primary, they must be finite. Thus, we have
\[ \rank_{\Op} (\ker (\al^{\vee})) - \rank_{\Op}(\coker (\al^{\vee}))= \rank_{\Op}\big(X_f(A/F_{\infty})_H\big) -\rank_{\Op}(X_f(A/F^{\cyc})).\]
 Combining this with the above calculations and applying Proposition \ref{Howson} for $X_f(A/F_{\infty})$, we obtain the required formula.
\epf

\section{Kida's formula} \label{Kida section}

We are in position to prove the main theorem of the paper.

\bt \label{Kida formula theorem}
Let $F_{\infty}$ be a strongly admissible pro-$p$ Lie extension of $F$ and $L$ a finite extension of $F$ contained in $F_{\infty}$ with $F^{\cyc}\cap L= F$. Let $\big(A, \{A_v\}_{v|p}, \{A^+_v\}_{v|\R} \big)$ be a datum satisfying $\textbf{(C1)-(C4)}$ and $\textbf{(Fin)}$. Assume that $X(A/F_{\infty})$ satisfies the $\M_H(G)$-property. Then we have
\[ \mu_{\Op\ps{\Ga_L}}\big(X(A/L^{\cyc})\big) = [L:F]\mu_{\Op\ps{\Ga}}\big(X(A/F^{\cyc})\big)\] and
\[ \rank_{\Op}(X_f(A/L^{\cyc})) = [L:F]\rank_{\Op}(X_f(A/F^{\cyc})) + \corank_{\Op}(A(L^{\cyc})) -[L:F] \corank_{\Op}(A(F^{\cyc}))  \] \[
+ \sum_{w\in R(L^{\cyc}/F^{\cyc})}\left([L:F]\corank_{\Op}\big(D_w(F^{\cyc}_w)\big) -\sum_{u|w}\corank_{\Op}\big(D_w(L^{\cyc}_{u})\big) \right).\]
Here $D_w$ denotes $A/A_w$ or $A$ accordingly as $w$ divides $p$ or not, and $R(L^{\cyc}/F^{\cyc})$ is set of primes of $F^{\cyc}$ consisting of all the primes above $p$ and all the ramified primes in $L^{\cyc}/F^{\cyc}$.
\et

\bpf
For a $\Op\ps{G}$-module $M$, one has the observation that $\mu_{\Op\ps{G}}(M) = [G:G']\mu_{\Op\ps{G'}}(M)$. Combining this with the first equality in Proposition \ref{rank formula}, we obtain the first equality of the theorem. On the other hand, one also has
$\rank_{\Op\ps{H}}(M) = [H:H_L]\rank_{\Op\ps{H_L}}(M)$. Applying this to $X(A/F_{\infty})$ and combining  with Proposition \ref{rank formula}, we have
\[ \rank_{\Op}(X_f(A/L^{\cyc})) - \corank_{\Op}(A(L^{\cyc})) + \sum_{\substack{u\in S(L^{\cyc}), \\ \dim H_{L,u}\geq 1}}\corank_{\Op}\big(D_u(L^{\cyc}_u)\big) = \]
\[ [L:F]\left(\rank_{\Op}(X_f(A/F^{\cyc})) - \corank_{\Op}(A(F^{\cyc})) + \sum_{\substack{w\in S(F^{\cyc}), \\ \dim H_w\geq 1}}\corank_{\Op}\big(D_w(F^{\cyc}_w)\big)\right).\]
Now since $L^{\cyc}/F^{\cyc}$ is a finite extension, it follows that for each $w\in S(F^{\cyc})$ and $u\in S(L^{\cyc})$ which lies above $w$, we have $\dim H_{L,u} \geq 1$ if and only $\dim H_{w} \geq 1$. Hence the summation $\displaystyle\sum_{\substack{u\in S(L^{\cyc}) \\ \dim H_{L,u}\geq 1}}$ can be rewritten as $\displaystyle\sum_{\substack{w\in S(F^{\cyc})\\ \dim H_{w}\geq 1}}\sum_{u|w}$. Now, upon rearranging the terms, we obtain
\[ \rank_{\Op}(X(A/L^{\cyc})) = [L:F]\rank_{\Op}(X(A/F^{\cyc})) + \corank_{\Op}(A(L^{\cyc})) -[L:F] \corank_{\Op}(A(F^{\cyc}))  \] \[
+ \sum_{w\in S(F^{\cyc})}\left([L:F]\corank_{\Op}\big(D_w(F^{\cyc}_w)\big) -\sum_{u|w}\corank_{\Op}\big(D_w(L^{\cyc}_{u})\big) \right).\]

Hence the required equality of this theorem will follow once we show that the summation $\displaystyle\sum_{w\in S(F^{\cyc})}$ can be rewritten as $\displaystyle\sum_{w\in R(L^{\cyc}/F^{\cyc})}$.
Clearly $R(L^{\cyc}/F^{\cyc})\subseteq S(F^{\cyc})$. Thus, it suffices to show that the term in the summand is zero when $w\in S(F^{\cyc})\setminus R(L^{\cyc}/F^{\cyc})$. But since such a prime does not divide $p$ and is not ramified in $L^{\cyc}/F^{\cyc}$, it must therefore split completely in $L^{\cyc}/F^{\cyc}$. Hence for each such prime $w$, we have $D_w(L^{\cyc}_u) =  D_w(F^{\cyc}_w)$ for every $u|w$, where there are $[L:F]$ of these $u$'s. Consequently, the term in the sum for this prime $w$ is zero. Thus, we have proven our result.
\epf

For the remainder of this section, we discuss how the above theorem can be used to recover Kida's formula for the $\mu=0$ case. Before doing so, we record the following preliminary observation.

\bl \label{burn-ven}
Let $F$ be a number field which contains a primitive $p$-th root of unity and $L$ a finite $p$-extension of $F$. Then there exists a strongly admissible pro-$p$ Lie extension $F_{\infty}$ of $F$ containing $L$.
\el

\bpf
This can be proven similarly to that in \cite[Lemma 6.1]{BV}.
\epf

We now apply Theorem \ref{Kida formula theorem} to recover Kida's formula for the $\mu=0$ case.

\bt \label{Kida formula theorem mu=0}
Let $F$ be a number field which contains a $p$-th root of unity and $L$ a finite $p$-extension of $F$. Let $\big(A, \{A_v\}_{v|p}, \{A^+_v\}_{v|\R} \big)$ be a datum satisfying $\textbf{(C1)-(C4)}$ and $\textbf{(Fin)}$. Suppose that $X(A/F^{\cyc})$ is finitely generated over $\Op$. Then $X(A/L^{\cyc})$ is finitely generated over $\Op$ and
\[ \rank_{\Op}(X(A/L^{\cyc})) = [L:F]\rank_{\Op}(X(A/F^{\cyc})) + \corank_{\Op}(A(L^{\cyc})) -[L:F] \corank_{\Op}(A(F^{\cyc}))  \] \[
+ \sum_{w\in R(L^{\cyc}/F^{\cyc})}\left([L:F]\corank_{\Op}\big(D_w(F^{\cyc}_w)\big) -\sum_{u|w}\corank_{\Op}\big(D_u(L^{\cyc}_{u})\big) \right),\]
where $R(L^{\cyc}/F^{\cyc})$ is set of primes of $F^{\cyc}$ consisting of primes above $p$ and the primes that are ramified in $L^{\cyc}/F^{\cyc}$.
\et

\bpf
Let $F_{\infty}$ be a strongly admissible pro-$p$ Lie extension $F_{\infty}$ of $F$ containing $L$ which is obtained via Lemma \ref{burn-ven}.  By a similar argument to that in \cite[Corollary 2.5]{PW06} or \cite[Corollary 3.4]{HM}, $X(A/L^{\cyc})$ is finitely generated over $\Op$. Also, it follows from \cite[Theorem 2.1]{CS12} that $X(A/F_{\infty})$ is finitely generated over $\Op\ps{H}$, where $H=\Gal(F_{\infty}/F^{\cyc})$. In particular, $X(A/F_{\infty})$ satisfies the $\M_H(G)$-property. Thus, Theorem \ref{Kida formula theorem} applies to yield the required conclusion.
\epf

\section{Review of elliptic curves} \label{elliptic review}

In this section, we revisit Kida's formula for the elliptic curve of Hachimori-Matsuno \cite{HM} and show that our theorem specializes to recover their theorem. For simplicity, we shall assume that $p\geq 5$. Let $E$ be an elliptic curve defined over a number field $F$ such that for every prime $v$ of $F$ above $p$, $E$ has either good ordinary reduction or split multiplicative reduction at $v$.  We shall also assume that $F$ contains all the $p$th roots of unity. Recall that by \cite[p.\ 150-151]{CG}, for each prime $v$ of $F$ above $p$, we have a short exact sequence
\[ 0\lra C_v\lra E(p) \lra D_v\lra 0 \]
of discrete $\Gal(\bar{F}_v/F_v)$-modules which is characterized by the facts that $C_v$ is divisible and $D_v$ is the maximal quotient of $E(p)$ by a divisible subgroup such that the inertia group acts on $D_v$ via a finite quotient. In fact, by our assumptions (also see \cite[pp. 69-70]{GR99}), $D_v$ has the following description
\[ D_v =\begin{cases}  \widetilde{E}_v(p),& \mbox{if $E$ has good ordinary reduction at $v$}, \\
      \Qp/\Zp, & \mbox{if $E$ has split multiplicative reduction at $v$}.\end{cases}\]
Here $\widetilde{E}_v$ is the reduced curve $E$ at $v$.

The following finiteness result of $\widetilde{E}_v(F^{\cyc}_w)(p)$ for a good ordinary prime $v$ is well-known but for a lack of reference, we include a proof here.

\bl \label{finiteness}
Let $E$ be an elliptic curve over $F$ with good ordinary reduction at the prime $v$ which lies above $p$. Then for every prime $w$ of $F^{\cyc}$ above $v$, we have that $\widetilde{E}_v(F^{\cyc}_w)(p)$ is finite.
\el

\bpf
The long cohomology exact sequence of
\[ 0\lra C_v\lra E(p) \lra \widetilde{E}_v(p)\lra 0 \]
 gives rise to an exact sequence
 \[ E(F^{\cyc}_w)(p)\lra \widetilde{E}_v(F^{\cyc}_w)(p)\lra H^1(F^{\cyc}_w, C_v)\lra H^1(F^{\cyc}_w, E(p)) \lra H^1(F^{\cyc}_w, \widetilde{E}_v(p)) \lra 0,\]
 where the final zero follows from that fact that $H^2(F^{\cyc}_w,C_v)=0$ (cf.\ \cite[Theorem 7.1.8(i)]{NSW}). Since our elliptic curve has good ordinary reduction at the prime $v$, Imai's theorem \cite{Imai} asserts that $E(F^{\cyc}_w)(p)$ is finite. On the other hand, a local Euler characteristics argument (cf.\ \cite[\S 3]{G89}) shows that $\corank_{\Zp}(H^1(F^{\cyc}_w, E(p))) =
 \corank_{\Zp}(H^1(F^{\cyc}_w, C_v)) + \corank_{\Zp}(H^1(F^{\cyc}_w, \widetilde{E}_v(p))$. Putting this information into the exact sequence, we see that the $\widetilde{E}_v(F^{\cyc}_w)(p)$ has trivial $\Zp$-corank.
\epf

Now for the local terms, we have
\[\corank_{\Zp}\big(D_w(F^{\cyc}_w)\big) = \begin{cases}  1,& \mbox{if $E$ has split multiplicative reduction at $w$}, \\
      2, & \mbox{if $w\nmid p$ and $E$ has good reduction at $w$ with $E(F^{\cyc}_w)[p]\neq 0$},
      \\
      0, & \mbox{otherwise,}\end{cases}
 \]
where we have made use of \cite[Proposition 5.1]{HM} for the primes not dividing $p$.

Let $L$ be a finite $p$-extension of $F$. By Lemma \ref{burn-ven}, one can find a strongly admissible pro-$p$ $p$-adic Lie extension $F_{\infty}$ of $F$ which contains $L$.
Denote by $P_0(F^{\cyc})$ the set of primes of $S(F^{\cyc})$, whose decomposition group of $H$ at $w$ has dimension $\geq 1$. Let $M(F^{\cyc})$ (resp., $P_1(F^{\cyc})$) be the set of split multiplicative primes of $E$ in $P_0(F^{\cyc})$ which lie above $p$ (resp., do not lie above $p$). Let $P_2(F^{\cyc})$ be the primes $w$ of $P_0(F^{\cyc})$ at which $w\nmid p$ and $E$ has good reduction at $w$ with $E(F^{\cyc}_w)[p]\neq 0$. Assuming that $X(E/F_{\infty})$ satisfies the $\M_H(G)$-conjecture, Theorem \ref{Kida formula theorem} then reads as
 \[ \rank_{\Zp}(X_f(E/L^{\cyc})) = [L:F]\rank_{\Zp}(X_f(E/F^{\cyc}))
+ \sum_{w\in M(F^{\cyc})\cup P_1(F^{\cyc})}\left([L:F] -\sum_{u|w}1\right) \]
\[+ 2\sum_{w\in P_2(F^{\cyc})}\left([L:F] -\sum_{u|w}1 \right), \]
where the $u$'s denote the primes in $S(L^{cyc})$, and noting that $E(F^{\cyc})(p)$ and $E(L^{\cyc})(p)$ are finite by a result of Ribet \cite{Ri}. One can see that the above agrees with that in \cite[Theorems 3.1 and 8.1]{HM} and \cite[Lemma 3.6]{HV}.

\br
One can also check that our Theorem \ref{Kida formula theorem} (or Theorem \ref{Kida formula theorem mu=0}) recovers \cite[Theorem 4.2]{CDLSS} and \cite[Theorem 2.8]{PW06}.
\er

We end by discussing an example taken from \cite{DD15}. Let $E$ be the elliptic curve defined by
\[ y^2 = x^3 -24z^7\sqrt{z+3} x^2 +zx,\]
where $z= \frac{\sqrt{5}-1}{2}$. This elliptic curve is defined over $\Q(\sqrt{z+3})$ and has good ordinary reduction at the prime above $5$. Let $F$ be the field $\Q(\zeta_{20})$ and $F_{\infty}$ the $\Z_5^5$-extension of $F$. It is shown in \cite{DD15} that the $\mu$-invariant of $X(E/F^{\cyc})$ is positive. Thus, if we assume that $\M_H(G)$-conjecture holds for $X(E/F_{\infty})$, then it follows from the above discussion that for every finite extension $L$ of $F$ contained in $F_{\infty}$, we have
\[ \rank_{\Zp}(X_f(E/L^{\cyc})) = [L:F]\rank_{\Zp}(X_f(E/F^{\cyc})).\]
Unfortunately, we do not have a way of establishing $\M_H(G)$-conjecture for $X(E/F_{\infty})$ as yet (but see \cite[Section 3]{CS12} and \cite[pp 1985]{DD15} for some partial results in this direction). It could be of interest to at least verify the above rank equality numerically (for certain $L$). We however do not know how to approach this problem at this point of writing.

\section{Remarks on the classical case}

We give a sketch how the discussion of Section \ref{Kida section} specializes to the classical case. From now on, $F$ will denote a totally real number field. Write $\Delta= \Gal(F(\mu_p)/F)$. Let
\[ \kappa: \Gal(F(\mu_{p^\infty})/F)\lra \mathrm{Aut}(\mu_{p^\infty})\cong \Zp^{\times} \]
be the cyclotomic character. Writing $d=|\Delta|$, we define
\[ e_{i} = \frac{1}{d}\sum_{\sigma\in \Delta} \kappa^{-i}(\sigma)\sigma^i \in \Zp[\Delta].\]
Let $i$ be an even integer such that $0\leq i\leq d$. Set $A=\Qp/\Zp(\kappa^{-i})$ and $A_v=\Qp/\Zp(\kappa^{-i})$ for all $v|p$, where
$\Qp/\Zp(\kappa^{-i})$ is $\Qp/\Zp$ as $\Zp$-module with a Galois action by $\sigma\cdot x = \kappa^{-1}(\sigma) x$. Then one can check that
\[S^{str}(A/F(\mu_{p^{\infty}})) = H^1\big(G_S(F(\mu_{p^{\infty}})),A\big) \cong \big(e_iG_S(F(\mu_{p^{\infty}}))^{ab}(p)\big)^{\vee}. \]
Here $G_S(F(\mu_{p^{\infty}}))^{ab}$ is the abelianization of the group $G_S(F(\mu_{p^{\infty}}))$. Note that $e_iG_S(F(\mu_{p^{\infty}}))^{ab}(p)$ is a torsion $\Zp\ps{\Ga}$-module by \cite[Proposition 11.4.5]{NSW}, where $\Ga=\Gal(F^\cyc/F)$.

We may proceed as in the previous section to obtain a Kida formula for $S^{str}(A/F(\mu_{p^{\infty}}))$ for every even $i$ under the assumption of the $\M_H(G)$-conjecture. Summing over these $i$'s, we obtain \cite[Corollary 11.4.11]{NSW}. Upon combining this with \cite[Corollary 11.4.4]{NSW}, we obtain the classical Kida's formula for the minus class groups \cite{Iw, Ki}. We should mention that Hachimori-Sharifi has proven the classical Kida's formula under a slight weakening of the usual $\mu=0$ assumption (cf. \cite[Theorem 2.1]{HS}). However, it would seem that their approach cannot be extended to prove the general case (see \cite[Remark 2.2.2]{Gr11}). Of course, we should mention that in this classical situation, Iwasawa has conjectured that $\mu=0$ always holds, although this has only been verified for abelian number fields (\cite{FW}).

\end{document}